\documentclass[12pt]{amsart}
\usepackage{amsmath,amssymb,amsfonts,epsfig}
\usepackage{xcolor}
\usepackage{amsfonts}
\usepackage{amssymb}
\usepackage{graphicx}
\usepackage{pstricks}
\usepackage{amsmath}
\usepackage{amscd}
\usepackage{epsfig}
\usepackage{bm}

\usepackage{amsxtra}
\usepackage{etoolbox}
\usepackage{mathabx}
\usepackage{hyperref}
\usepackage{cleveref}

\usepackage{tkz-tab}

\usepackage[all]{xy}

\usepackage{cancel}
\usepackage{ulem}
\usepackage{xcolor}

\newtheorem{thm}{Theorem}[section]
\newtheorem{prop}[thm]{Proposition \!\!}

\newtheorem{lem}[thm]{Lemma \!\!}

\newtheorem{definition}[thm]{Definition}

\newcommand\mk{\medskip}
\newcommand\bk{\bigskip}

\newcommand{\pfend}{\hfill $\Box$ \medskip}

\begin{document}


\title{Toeplitz operators on some function spaces}
\author{Chafiq Benhida}
\address{UFR de Math\'{e}matiques, Universit\'{e} des Sciences et
Technologies de Lille, F-59655, Villeneuve-d'Ascq Cedex, France}
\email{chafiq.benhida@univ-lille.fr}

\author{George R. Exner}
\address{Department of Mathematics, Bucknell University, Lewisburg, Pennsylvania 17837, USA}
\email{exner@bucknell.edu}

\author{Ji Eun Lee}
\address{Department of Mathematics and Statistics, Sejong University, Seoul 143-747, Korea}
\email{jieunlee7@sejong.ac.kr}

\author{Jongrak Lee}
\address{Department of Mathematics, Sungkyunkwan University, Suwon 16419, Korea}
\email{jrlee01@skku.edu}



\begin{abstract}
We reconsider studies of Toeplitz operators on function spaces (the weighted Bergman space, the generalized derivative Hardy space) and the H-Toeplitz operators on the Bergman space.  Past studies have considered the presence or absence of hyponormality or questions of contractivity or expansivity;  we provide structure theorems for these operators that allow us to recapture, and often considerably improve, these results.  In some cases these operators or their adjoints are actually in more restrictive classes, such as subnormal or moment infinitely divisible ($\mathcal{MID}$).\
\end{abstract}

\maketitle

\setcounter{tocdepth}{4}

\setcounter{secnumdepth}{4}

\tableofcontents

\setcounter{tocdepth}{3}

\medskip

\section{Introduction and Statement of Main Results} \label{Intro}

\medskip

The study of Toeplitz operators on the Hardy space or Bergman space is long-standing;   more recently Toeplitz operators on other function spaces have been considered.  Recall that the unifying features are a functional Hilbert space, a privileged subspace (often consisting of analytic functions), and multiplication by a function -- the symbol -- followed by a projection onto the privileged subspace.  In particular, with all definitions to be reviewed below, spaces defined on the unit disk $\mathbb{D}$ such as the weighted Bergman spaces (sometimes called the $\alpha$-Bergman spaces or weighted Hardy spaces; see \cite{HKZ}, \cite{SkL1}, \cite{LS}), the derivative Hardy space or generalized derivative Hardy spaces (\cite{GL}, \cite{KLL}), or the so-called H-Toeplitz operators on the Bergman space obtained by use of the space of harmonic functions on $\mathbb{D}$ (see, e.g., \cite{GS} and \cite{SkL}) have been explored. Often the questions considered have been the hyponormality, or its absence, of the operator and the conditions yielding contractive or expansive operators.

Using a reasonably unified approach that represents, for some symbols,  these operators (or perhaps their adjoints) as weighted shifts, or direct sums including weighted shifts, we provide structure theorems that, in many cases, improve hyponormality to subnormality or the even better property of $\mathcal{MID}$, render more transparent hyponormality or its absence, and similarly display the conditions for contractivity or expansivity more directly.  Results from \cite{BCE1} and especially \cite{BCE2} play an important role in this approach.

\vspace{.1in}

We turn next to setting notation and recalling definitions.

\medskip

\subsection{The function spaces}\hfill\\

 Let $\mathbb{D}$ denote the open unit disk in the
complex plane $\mathbb{C}$ and $dA$ the area measure on $\mathbb{C}$. Let $\mathbb{R}_+$ be the set of positive real numbers and $\mathbb{N}$ be the set of positive integers;  we set $\mathbb{N}_n = \{j \in \mathbb{N}: j \geq n \ \text{for} \ n \in \mathbb{N}\}$, and $\mathbb{N}_0 = \mathbb{N} \cup \{0\}$.  The space
$L^2(\mathbb{D})$ is a Hilbert space with the inner product
$$
\langle f, \ g \rangle = \frac{1}{\pi} \int_{\mathbb{D}} f(z)
\overline{g(z)} dA(z).
$$
The Bergman space $L^2_a(\mathbb{D})$ consists of all analytic functions in $L^2(\mathbb{D})$ and $L^{\infty}(\mathbb{D})$ is the space of the essentially bounded measurable function on
$\mathbb{D}$. For $\varphi \in L^{\infty}(\mathbb{D}),$ the multiplication
operator $M_{\varphi}$ on $L^2(\mathbb{D})$ is defined by
$M_{\varphi}(f)=\varphi \cdot f$ and the  Toeplitz operator $T_{\varphi}$ on $L^2_a(\mathbb{D})$ is defined by
$$
T_{\varphi} (f) = P({\varphi}\cdot f),
$$
where $P$ denotes the orthogonal projection of $L^2 (\mathbb{D})$ onto $L^2_a(\mathbb{D})$ and $f\in L^2_a(\mathbb{D})$. It is clear
that those operators are bounded if $\varphi \in L^{\infty}(\mathbb{D}).$

The harmonic Bergman space $L^2_{\it{harm}}(\mathbb{D})$  denotes the space of all complex-valued harmonic functions in $L^2(\mathbb{D})$.  The
space $L^2_{\it{harm}}(\mathbb D)$ is a closed subspace of $L^2(\mathbb{D})$ and it is a Hilbert space.
Let $P_{\it{harm}}$ be the orthogonal projection from the space $L^2(\mathbb{D})$
onto the space $L^2_{\it{harm}}(\mathbb{D})$.

\medskip

\subsubsection{Weighted Bergman spaces} \hfill \\

For $-1 <\alpha<
\infty$, the {\it weighted Bergman space} $A_{\alpha}^2(\mathbb D)$ is the
space of all analytic functions in $L^2(\mathbb D, dA_{\alpha})$, where
$$
dA_{\alpha}(z)=(\alpha+1)(1-|z|^2)^{\alpha}dA(z).
$$
If $f, \ g \in L^2(\mathbb D, dA_{\alpha})$, we write
$$
\langle f, \ g \rangle=\int_{\mathbb
	D}f(z)\overline{g(z)}dA_{\alpha}(z) \quad \text{and} \quad \|f\|=\left(\int_{\mathbb
	D}|f(z)|^2dA_{\alpha}(z)\right)^{\frac{1}{2}}.
$$
The space $L^2(\mathbb D, dA_{\alpha})$ is a Hilbert space
with the above inner product.
For a nonnegative integer $n$ and any $z \in \mathbb D$, let
$$
e_n(z)=\sqrt{\frac{\Gamma(n+\alpha+2)}{\Gamma(n+1)\Gamma(\alpha+2)}}~z^n,
$$
where $\Gamma(\cdot)$ is the usual Gamma function. Then $\{e_n\}$ is an orthonormal basis for
the weighted Bergman space (cf. \cite{HKZ}).
Given a bounded measurable function $\varphi
\in L^{\infty}(\mathbb D),$ the \textbf{Toeplitz operator} $T_{\varphi}$ with symbol
$\varphi$ on $A_{\alpha}^2(\mathbb D)$
is defined by
$$
T_{\varphi}g := P({\varphi}\cdot g) \ \ (g\in A_\alpha^2(\mathbb D))
$$
\noindent where $P$ denotes the usual orthogonal projection from $L^2
(\mathbb D, dA_\alpha)$ onto $A_\alpha^2(\mathbb D)$.

We record a result of $P$ on $A^2_{\alpha}(\mathbb D)$ crucial for the sequel.

\begin{lem}{\label{Lemma 1}} {\rm(\cite{HL2})}
 For any nonnegative integers $s, t,$
		$$
		P(\overline{z}^t z^s)=\begin{cases}
			\frac{\Gamma(s+1)\Gamma(s-t+\alpha+2)}{\Gamma(s+\alpha+2)\Gamma(s-t+1)}
			z^{s-t} \qquad \,\, \mbox{\rm if} \ s \geq t \\ 0 \qquad \qquad \qquad \qquad
			\quad \quad  \mbox{\rm if} \ s < t.
		\end{cases}
		$$
\end{lem}

\medskip

\subsubsection{H-Toeplitz operators on the Bergman space}  \hfill \\

We consider an alternative to the classical Toeplitz operators on the Bergman space.   Let $s,t$ be nonnegative integers and $P$ be the orthogonal projection from  $L^2(\mathbb D)$ to  $L^2_a(\mathbb D)$. Then we have
$$
P(\overline{z}^t z^s)=\begin{cases} \frac{s-t+1}{s+1} z^{s-t} \qquad
	\text{if} \ s \geq t \\ 0 \qquad \qquad  \hspace{.255in} \text{if} \ s < t_.
\end{cases}
$$
\medskip
The following result concerning $P_{harm}$ will be used frequently.
\begin{lem}\label{lem2.1}{\rm(\cite{GS})}
	In the harmonic Bergman space $L^2_{\it{harm}}(\mathbb D)$, for nonnegative integers $s$ and $t$, we have
\begin{equation}  \label{eg:defPharm}
P_{\it{harm}}(\overline{z}^t z^s)=\begin{cases} \frac{s-t+1}{s+1} z^{s-t} \qquad \text{{\rm if}} \ s \ge t \\
		\frac{t-s+1}{t+1} \overline{z}^{t-s} \qquad \text{{\rm if}} \ s < t.
	\end{cases}
\end{equation}

\end{lem}

In order to define the notion of an H-Toeplitz operator on $L^2_a(\mathbb D)$, we first consider the operator $K : L^2_a(\mathbb D) \rightarrow L^2_{\it{harm}}(\mathbb D)$ defined by
\begin{equation*}  \label{eq:defK}
K(e_{2n}(z)) = e_n(z) = \sqrt{n+1} \, z^n
\ \text{and} \   K(e_{2n+1}(z)) = \overline{e_{n+1}(z)} = \sqrt{n+2} \, \, \overline{z}^{n+1}
\end{equation*}
for all $n\ge 0$ and $z\in \mathbb D$. It can be checked that the operator $K$ is bounded
and linear on $L^2_a(\mathbb D)$ with $\|K\| = 1$. Further, the adjoint $K^*$ of the operator $K$ is
given by
$$
K^*(e_{n}(z)) = e_{2n}(z)
\  \, \, \text{and} \   K^*(\overline{e_{n+1}(z)}) = e_{2n+1}(z) \quad
$$
for all  $n\ge 0$. From the definition of $K$ and $K^*$, we have that $KK^*=I_{L^2_{\it{harm}}(\mathbb D)}$ and $K^*K=I_{L^2_a(\mathbb D)}$.  Using the definitions of $K$ and $K^*$, we have that
\begin{eqnarray}
	K(z^{2n})&=& \frac{\sqrt{n+1}}{\sqrt{2n+1}}z^n, \nonumber \\
    K(z^{2n+1})&=&\frac{\sqrt{n+2}}{\sqrt{2n+2}}\overline{z}^{n+1}, \nonumber \\
    K^*(z^{n})&=&\frac{\sqrt{2n+1}}{\sqrt{n+1}}z^{2n}, \, \mbox{\rm and} \label{eq:Kstarzupn} \\
	K^*(\overline{z}^{n})&=&\frac{\sqrt{2n}}{\sqrt{n+1}}{z}^{2n-1}.  \label{eq:Kstarzbarupn}
\end{eqnarray}

\medskip

We may now define H-Toeplitz operators on the Bergman space $L^2_a(\mathbb D)$ using the operator $K$.
\begin{definition}{\rm(\cite{GS})}
	For $\varphi\in L^{\infty}(\mathbb D)$, the H-Toeplitz operator $B_{\varphi}$ with the symbol $\varphi$ is defined as the operator
	$B_{\varphi} : L^2_a(\mathbb D) \rightarrow L^2_a(\mathbb D)$ such that $B_{\varphi}(f) = PM_{\varphi}K(f)$ for all $f \in L^2_a(\mathbb D)$.
\end{definition}

\medskip

The next proposition follows from the definition of the H-Toeplitz operators.
\begin{prop}{\rm(\cite{GS})}
	For $\varphi, \psi \in L^{\infty}(\mathbb D)$, the operator $B_{\varphi}$ satisfies the following:
	
	{\rm(i)} $B_{\varphi}$ is a bounded linear operator on $L^2_a(\mathbb D)$ with $\|B_{\varphi}\|\le \|\varphi\|_{\infty}$.
	
	{\rm(ii)} For any scalar $\alpha$ and $\beta$, $B_{\alpha \varphi+\beta\psi}=\alpha B_{\varphi}+\beta B_{\psi}$.
	
	{\rm(iii)} The adjoint of the H-Toeplitz operator $B_{\varphi}$ is given by $B^*_{\varphi}=K^*P_{\it{harm}}M_{\overline{\varphi}}$.
\end{prop}

\medskip

\subsubsection{Generalized derivative Hardy space} \hfill \\

In 2018, C. Gu and S. Luo (\cite{GL}) introduced the derivative Hardy space $S^2_1(\mathbb D)$ as follows:\\
\begin{eqnarray*}
S^2_{1}(\mathbb{D})&=& \left\{f\in H(\mathbb{D}) : \|f\|^2_{S^2_{1}}=
|f\|^2_{H^2}+{\frac{3}{2}}\|f'\|^2_{A^2}+\frac{1}{2}\|f'\|^2_{H^2}<\infty\right\}\\
&=&\left\{f\in H(\mathbb{D}) : \|f\|^2_{S^2_{1}}=\sum_{n=0}^{\infty}\frac{(n+1)(n+2)}{2}|f_n|^2<\infty\right\}.
\end{eqnarray*}

The reproducing kernel of the derivative Hardy space  $S^2_1(\mathbb D)$ is $$K_w(z)=\frac{2}{(\overline{w}z)^2}\left(\overline{w}z+(\overline{w}z-1)\ln \frac{1}{1-\overline{w}z}\right).$$

Recently, E. Ko, J. E. Lee and J. Lee (\cite{KLL}) defined  the generalized derivative Hardy space $S^2_{\alpha,\beta}(\mathbb D)$ { for $\alpha,\beta\in{\mathbb N}$ with {$\alpha<\beta$}} as

\begin{eqnarray*}
S^2_{\alpha,\beta}(\mathbb{D})&=&\left\{f\in H(\mathbb{D}) : \|f\|^2_{S^2_{\alpha,\beta}}=\|f\|^2_{H^2}+{\frac{{\alpha+\beta}}{\alpha\beta}}\|f'\|^2_{A^2}+\frac{1}{\alpha\beta}\|f'\|^2_{H^2}<\infty\right\}\\
&=&\left\{f\in H(\mathbb{D}) : \|f\|^2_{S^2_{\alpha,\beta}}=\sum_{n=0}^{\infty}\frac{(n+\alpha)(n+\beta)}{\alpha\beta}|f_n|^2<\infty\right\}.
\end{eqnarray*}

\medskip

We will need the following lemma.

\smallskip

\begin{lem}{\rm (\cite{KLL})}    \label{eq:PforgenderivHS}
  For any nonnegative integers $s, t$, we have
\begin{eqnarray*}
&& {(a)}~ \langle{z}^t, z^s\rangle=\begin{cases} \frac{(s+\alpha)(s+\beta)}{\alpha\beta}
 \qquad  & \text{if} \ s = t \\ 0 \quad \quad  \, \ \
& \text{if} \ s \neq t,~\mbox{and}
\end{cases}\\
&&{(b)}~ P(\overline{z}^t z^s)=\begin{cases} \frac{(s+\alpha)(s+\beta)}{(s-t+\alpha)(s-t+\beta)} z^{s-t} \qquad
&\text{if} \ s \geq t \\ 0 \qquad \qquad  \quad\, \,
 \, \, &\text{if} \ s < t.
\end{cases}
\end{eqnarray*}
\end{lem}

One computes that an orthonormal basis is given by $\{e_n\}_{n=0}^\infty$ where
\begin{equation}   \label{eq:esubnforgenderivHS}
e_n = \sqrt{\frac{\alpha \beta}{(n + \alpha)(n + \beta)}} z^n.
\end{equation}

\medskip

\subsection{Unilateral weighted shifts, their properties, and classes}  \label{subse:shiftspropertiesclasses}\hfill\\

We employ the standard notation for weighted shifts: let  $\ell^2$ be the classical Hilbert space of square summable complex sequences, with canonical orthonormal basis $e_0, e_1, e_2,  \ldots$ (note indexing begins at zero). \  Let $\alpha: \alpha_0, \alpha_1, \alpha_2, \ldots$ be a (bounded) positive \textbf{weight sequence} and  $W_\alpha$ be  the weighted shift defined by $W_\alpha e_j := \alpha_j e_{j+1} \;\; (j =0,1,2, \ldots)$ and extended by linearity. \ (While weighted shifts can be defined for any bounded sequence $\alpha$, without loss of generality for our questions of interest we can and do assume henceforth that $\alpha$ is positive.) \  \ The \textbf{moments} $\gamma = (\gamma_n)_{n=0}^\infty$ of the shift are defined by $\gamma_0 := 1$ and $\gamma_n := \prod_{j=0}^{n-1} \alpha_j^2$ for $n \geq 1$. \ It is well known from \cite[III.8.16]{Con} and \cite{GW} that a weighted shift $W_\alpha$ is subnormal if and only if it has a \textbf{Berger measure}, meaning a probability measure $\mu$ supported on $[0, \|W_\alpha\|^2]$ such that
$$
\gamma_n = \int_0^{\|W_\alpha\|^2} t^n d \mu(t), \hspace{.2in} n = 0, 1,2, \ldots .
$$

Recall that a bounded operator $T$ is hyponormal if $T^*T-TT^* \ge 0$, and one computes easily that for a weighted shift $W_{\alpha}$ this is exactly $\alpha_n^2 \le \alpha_{n+1}^2$ for all $n =0,1,\ldots$.  An operator is $k$-hyponormal (some $k = 1, 2, \ldots$) if a certain $(k+1) \times (k+1)$ operator matrix is positive semi-definite (see, for example, \cite{Br}).
 It is well known from \cite[Theorem 4]{Cu} that for weighted shifts $k$--hyponormality reduces to the positivity, for each $n$, of the $(k+1) \times (k+1)$ (Hankel) moment matrix $M_{\gamma}(n,k)$, where
$$
M_{\gamma}(n,k) = \left(
\begin{array}{cccc}
\gamma _{n} & \gamma _{n+1} & \cdots & \gamma _{n+k} \\
\gamma _{n+1} & \gamma _{n+2} & \cdots & \gamma _{n+k+1} \\
\vdots & \vdots & \ddots & \vdots \\
\gamma _{n+k} & \gamma _{n+k+1} & \cdots & \gamma _{n+2k}%
\end{array}
\right).
$$

\medskip

We pause briefly to give some sequence and function monotonicity conditions useful for weighted shifts.  Given a sequence $a = (a_j)_{j=0}^\infty$, let $\nabla$ (the \textbf{forward difference operator}) be defined by $(\nabla a)_j := a_j - a_{j+1}$ for all $j$, and the iterated forward difference operators $\nabla^{n}$ by $
\nabla^{0} a := a \; \; \textrm{ and } \; \; \nabla^{n} := \nabla (\nabla^{n-1}),$
for $n \geq 1$. \ For example, $
(\nabla^2a)_j=a_j-2a_{j+1}+a_{j+2} \; \; (j =0,1,\ldots)$.

A sequence is $n$\textbf{--alternating} if $(\nabla^{n} a)_k \leq 0$ for all $k = 0, 1,2, \ldots$, and \textbf{completely alternating} if it is $n$--alternating for all $n = 1, 2, \ldots$. \  A sequence is \textbf{log completely alternating} if the sequence $(\ln a_j)$ is  completely alternating.  It is shown in \cite{BCE1} that if a sequence is completely alternating, it is log completely alternating, and that the reverse is not true in general.    A function $f : \mathbb{R}_+ \rightarrow \mathbb{R}_+$ is a \textbf{Bernstein} function if $f' \geq 0$ and $f^{(2n)} \leq 0$ ($n = 1, 2, \ldots$) and $f^{(2n - 1)} \geq 0$, ($n = 1, 2, \ldots$).  We define $f$ to be  \textbf{log Bernstein} if $\ln f$ is  Bernstein but with the relaxation that we do not require $\ln f$ to be nonnegative, but insist only on the appropriate alternation of signs of the derivative. \  It is well known that if a sequence is interpolated by a Bernstein (respectively, log Bernstein) function, then the sequence is completely alternating (respectively, log completely alternating).  Finally, a Hilbert space operator is \textbf{completely hyperexpansive} (see \cite{At}) if a certain operator inequality holds;  for a weighted shift, this simplifies to
$$\sum_{i=0}^n (-1)^i \binom{n}{i} \gamma_{i + j} \leq 0, \quad j\in \mathbb{N}_0, n \in \mathbb{N}.$$
An operator is \textbf{$m$-alternatingly hyperexpansive} (see \cite{SA}) if a different operator inequality holds;  for shifts, this becomes
\begin{equation*} \label{eq:defmAltHE}
\sum_{i=0}^m (-1)^{(m-i)} \binom{m}{i} \gamma_{i + j} \geq 0, \quad j \in \mathbb{N}_0;
\end{equation*}
a shift is \textbf{alternatingly hyperexpansive} if it is $m$-alternatingly hyperexpansive for all $m \in \mathbb{N}$.

Returning to the considerations of shifts, in \cite{BCE1} and \cite{BCE3} the authors considered a class of ``better than subnormal'' weighted shifts. \  We say that a shift $W_\alpha$ is \textbf{moment infinitely divisible} ($\mathcal{MID}$) if, for every $s > 0$, the ``Schur $s$--th power'' shift with weight sequence $(\alpha_n^{s})$ is subnormal. \ (Throughout this paper, we use ``Schur product,'' in any setting -- matricial, sequential, \ldots -- to mean the entrywise product of the objects at hand.)   It is easily seen that raising weights to the power $s$ and raising moments to the power $s$ are equivalent. \  Examples of $\mathcal{MID}$ shifts include the Agler shifts $A_j$, $j = 1, 2, \ldots$ (where $A_j$ has the weight sequence $\sqrt{\frac{n+1}{n+j}}$), with $j=1$ yielding the unweighted shift and $j=2$ the familiar Bergman shift \cite{BCE1}. \  This can be generalized to the fact that the so-called homographic shifts $S(a,b,c,d)$ (see \cite{CD}) with weight sequence $\sqrt{\frac{a n + b}{c n + d}}$ where $a$, $b$, $c$, and $d$ are positive real numbers such that $ad - bc > 0$, are in fact $\mathcal{MID}$  (see first \cite{CPY} and \cite{BCE1}).  The $\mathcal{MID}$ shifts are those with weights squared log completely alternating or, equivalently, with moments log completely monotone (\cite{BCE1} and \cite{BCE3} respectively).

Recently there has been study of a particular class of weighted shifts with representatives in various of the classes related to subnormality.
\begin{definition}[\cite{BCE4}]
For $p > 0$, $-1 < N < 1$ and $-1 < D < 1$, we define the \textbf{geometrically regular weighted shift} with parameters $p$ and $(N,D)$ (often \textbf{GRWS}) to be the weighted shift with weight sequence $\alpha = (\alpha_n)_{n=0}^\infty$, where
$$\alpha_n = \sqrt{\frac{p^n + N}{p^n + D}}, \quad n = 0, 1,2, \ldots\, .$$
\end{definition}

Fixing $p$ for the moment, \cite{BCE4} considered the GRWS for the parameter pair $(N,D)$ in the ``magic square'' $-1 < N < 1$ and $-1 < D < 1$.  Various examples of these shifts will appear in Subsection \ref{subse:HToepinBerg}, and the following theorem captures their representation in various classes, using a diagram of the magic square in Figure 1.

\bk

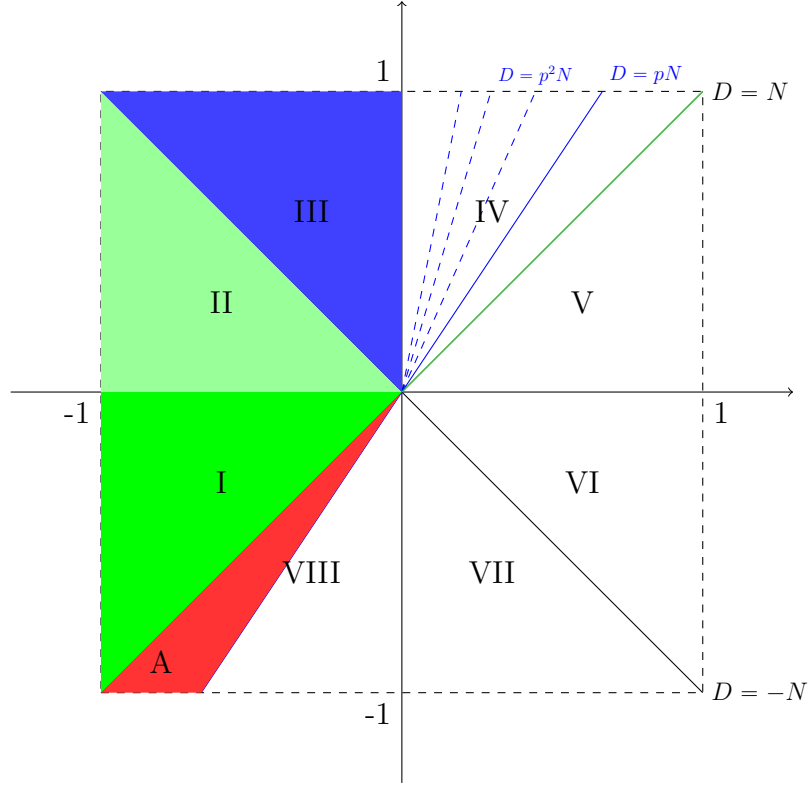
\begin{figure}


\begin{tikzpicture}[scale=4]

\draw[->] (-1.3,0) -- (1.3,0) ;

\draw[->] (0,-1.3) -- (0,1.3) ;

\draw  [-, dashed] (-1,-1)--(1,-1) --(1,1)--(-1,1)--cycle  ;

\draw (1,0) node[below right] {1} ;

\draw (-1,0) node[below left] {-1} ;

\draw (0,1) node[above left ] {1} ;

\draw (0,-1) node[below left] {-1} ;


\draw [][  domain=-1:1] plot(\x,{\x })node[right, scale=0.8] {$D=N$} ;

\draw [][  domain=-1:1] plot(\x,{-\x })node[right, scale=0.8] {$D=-N$} ;

\draw [blue ][  domain=-2/3:2/3] plot(\x,{3/2*\x })node[above right, scale=0.65] {$D=pN$} ;

\draw [-,blue, dashed][  domain=0:4/9] plot(\x,{9/4*\x })node[above , scale=0.6] {$D=p^2N$} ;

\draw [-, blue, dashed][  domain=0:8/27] plot(\x,{27/8*\x })node[above, scale=0.8] {} ;

\draw [-, blue, dashed][  domain=0:16/81] plot(\x,{81/16*\x })node[above, scale=0.8] {} ;

\fill[color=red!80] (0,0) -- (-1,-1) -- (-2/3,-1)  -- cycle ;

\fill[color=green] (0,0) -- (-1,-1) -- (-1,0)  -- cycle ;

\fill[green!40] (0,0) -- (-1,0) -- (-1,1)  -- cycle ;

\fill[blue!75] (0,0) -- (0,1) -- (-1,1)  -- cycle ;

\draw [green] (0,0) -- (1,1);

\draw (-0.6, -0.3) node [ ] {I};

\draw (-0.6, 0.3) node [ ] {II};

\draw (-0.3, 0.6) node [ ] {III};

\draw (0.3, 0.6) node [ ] {IV};

\draw (0.6, 0.3) node [ ] {V};

\draw (0.6, -0.3) node [ ] {VI};

\draw (0.3, -0.6) node [ ] {VII};

\draw (-0.3, -0.6) node [ ] {VIII};

\draw (-0.8, -0.9) node [ ] {A};


\end{tikzpicture}

\caption{Magic Square}
\label{MagicSquare}
\end{figure}

\begin{thm}  \label{th:detailsofmagicsquare}
Let $p > 1$ and $(N, D)$ such that $-1 < N < 1$ and $-1 < D < 1$. \  Let $W_\alpha$ be a geometrically regular weighted shift with parameters $p$ and $(N,D)$; that is,
$$\alpha_n :=  \sqrt{\frac{p^n + N}{p^n + D}}, \quad n = 0, 1,2, \ldots .$$
Then $W_\alpha$ has the following properties in the sectors given in the diagram:
\begin{enumerate}
      \item Along the main diagonal $D=N$, $W_\alpha$ is the (unweighted) unilateral shift and thereby both $\mathcal{MID}$ and completely hyperexpansive;
     \mk \item In Sector {\rm I}, the weight squared sequence $\alpha$ is interpolated by a Bernstein function and $W_\alpha$ is $\mathcal{MID}$;
      \mk \item In Sector {\rm II}, the weight squared sequence $\alpha$ is interpolated by a log Bernstein function and $W_\alpha$ is $\mathcal{MID}$;
      \mk \item In Sector {\rm III}, $W_\alpha$ is subnormal, and at least one shift in this sector is not $\mathcal{MID}$;
      \mk \item In Sector {\rm IV}, along the rays $D = p^k N$, $k = 0, 1, 2, \ldots$, $W_\alpha$ is subnormal with finitely atomic Berger measure and is not $\mathcal{MID}$ except for $N=D$;
      \mk \item In Sector {\rm IV}, in the subsector $p^{k-1}N < D < p^k N$, $W_\alpha$ is $k$-hyponormal but not $(k+1)$-hyponormal;
      \mk \item In Sector {\rm VIIIA}, $W_\alpha$ is completely hyperexpansive.
    \end{enumerate}
\end{thm}

We must make some additions to the results there since we will need to consider shifts of this form but whose parameters depart from the magic square.  First, if we consider the space ``above Sector III'' (so consisting of $(N, D)$ for which $D \geq 1$ and $-1 < N \leq 0$), the shifts arising from these parameters are subnormal.  The reason is that the form of the measure obtained for the shifts in Sector III works equally well here -- see \cite{BCE4}, and there is no problem with convergence of the relevant series for $N$ in the given range.  Similarly, subnormality of the shifts arising from points on the extensions of the ``special lines'' $D = p^k N$ (some $k = 0, 1, 2, \ldots$) outside the magic square are again subnormal;  the measure candidate used in \cite{BCE4} works equally here, and there are no problems with convergence (for $N \geq 1$) because the measure is finitely atomic.  Finally, shifts for points in the subsectors of Sector IV bounded by the special lines, but outside the magic square, have the same ``$k$- but not $(k+1)$-hyponormality'' property, because this argument is based upon determinant calculations (see the proof of Lemma 2.15 in \cite{BCE4}) which are unaffected by departure from the magic square.

\medskip

\subsection{Statement of main results}\hfill\\

\medskip

 In what follows we give somewhat abbreviated versions of the principal results shorn of certain technicalities;  the results in Section \ref{se:2} contain the complete statements.

For H-Toeplitz operators we have the following.

\begin{thm}  \label{briefth:seqt}
Let $B_{z^s \overline{z}^s}$ be an H-Toeplitz operator on the Bergman space with $s > 0$.  Then if $s =2^k - 1$ for some $k \in \mathbb{N}$, $B^*_{z^s \overline{z}^s}$ is  the direct sum of a one-dimensional normal operator and an infinite family of weighted shifts each subnormal.  For specified ranges in $s$, $B^*_{z^s \overline{z}^s}$ is the direct sum of a one-dimensional normal operator and an infinite family of weighted shifts each at least $k$-hyponormal.    In either case $B_{z^s \overline{z}^s}$ is not hyponormal;  both $B^*_{z^s \overline{z}^s}$ and $B_{z^s \overline{z}^s}$ are contractive.
\end{thm}

Remark that a result follows for certain H-Toeplitz operators with symbol a sum of those above:  see Theorem \ref{th:sumswhendeq0}.

\begin{thm}
Consider some H-Toeplitz operator $B_{z^t \overline{z}^s}$ on the Bergman space with $s > 0$ and $d = s-t > 0$.  If $s \leq 3d-2$ then $B^*_{z^t \overline{z}^s}$ is  a direct sum of weighted shifts each of which is $\mathcal{MID}$.  For any $s > 0$ and $d = s-t > 0$, $B^*_{z^t \overline{z}^s}$ is a direct sum of weighted shifts each subnormal, and is therefore subnormal.  It follows that in either case $B_{z^t \overline{z}^s}$ is not hyponormal.
\end{thm}

Again a result follows for certain sums: see Theorem \ref{th:sumsofHNliked}.

\begin{thm}
Consider $B_{z^t \overline{z}^s}$, a H-Toeplitz operator on the Bergman space with $t > s$ and set $\delta = t-s$.  If $\delta > 1$, then $B^*_{z^t \overline{z}^s}$ has a direct sum decomposition as the sum of an operator $M$ on $\mathbb{C}^{2 \delta}$ neither hyponormal nor co-hyponormal, a normal operator on $\mathbb{C}$, and a direct sum of weighted shifts each at least $1$-hyponormal, and at least $k$-hyponormal for certain ranges in the parameter $s$. Thus neither $B^*_{z^t \overline{z}^s}$ nor $B_{z^t \overline{z}^s}$ is hyponormal;  each is contractive.

In the case $\delta = 1$, there is the same form of decomposition, the operator $M$ on $\mathbb{C}^{2}$is normal, and $B^*_{z^t \overline{z}^s}$ is hyponormal;  $B_{z^t \overline{z}^s}$ is not hyponormal, and each is contractive.
\end{thm}

For Toeplitz operators on the weighted Bergman spaces we obtain the following, where again a result concerning sums follows as well (see Theorem \ref{thm:sumsforwtdBerg}).

\begin{thm}
Let $T_{\overline{z}^sz^{s+d}}$ be a Toeplitz operator on one of the  weighted Bergman spaces with $s \geq 0$ and $d \geq 0$.  If $d = 0$ the operator is diagonal and normal.  If $d > 1$, the operator is a direct sum of $d$ $\mathcal{MID}$ and hence subnormal weighted shifts. The operator is contractive, and except in the case $d = 0$, its adjoint is not hyponormal.
\end{thm}

Finally, for certain Toeplitz operators on the generalized derivative Hardy space (which are already known not to be even hyponormal unless normal  -- see \cite[Theorem 3.2]{KLL}) we obtain a shift-sum decomposition akin to those above with information about the included shifts.  In a certain special case these are Schur products of shifts which are $n$-alternatingly hyperexpansive in the sense of \cite{SA}.

\medskip

\section{Proofs and Discussion}  \label{se:2}

We consider in succession the H-Toeplitz operators, the Toeplitz operators on the weighted Bergman spaces, and those on the generalized derivative Hardy spaces.

\subsection{H-Toeplitz operators on the Bergman space}  \label{subse:HToepinBerg}\hfill\\

\medskip

We begin with an illustrative example to indicate the sort of structure result to be obtained, where the relevant calculations are simpler than in the abstract case.

\subsubsection{A motivating example}  \hfill \\

Suppose we consider the symbol $\varphi(z) = z^2 \overline{z}$, the multiplication operator $M_{z^2 \overline{z}}$ and the resulting H-Toeplitz operator $B_{z^2 \overline{z}} = PM_{z^2 \overline{z}}K$;  we will in fact concentrate on the adjoint $B^*_{z^2 \overline{z}} =  K^* P_{harm} M_{\overline{z}^2 z}$.  Using the results in \eqref{eg:defPharm} and \eqref{eq:Kstarzbarupn} we obtain that
\begin{eqnarray*}
B^*_{z^2 \overline{z}}(1) &=& B^*_{z^2 \overline{z}}(z^0)\\
&=&K^* P_{harm} M_{\overline{z}^2 z}(z^0) \\
&=& K^* P_{harm}(\overline{z}^2 z) \\
&=& K^*\left(\frac{2 - 1 + 1}{2+1} \overline{z}\right) \\
&=& \frac{2}{3} K^*(\overline{z}) \\
&=& \frac{2}{3} z.
\end{eqnarray*}
Recalling the relationship between the orthonormal basis $(e_n)_{n=0}^\infty$ and the $(z^n)_{n=0}^\infty$, namely $e_n = \sqrt{n+1} \,  z^n$ (see, e.g., \cite{HKZ}),  under the action of $B^*_{z^2 \overline{z}}$ we have $e_0 \mapsto \frac{2}{3}\cdot \frac{1}{\sqrt{2}} e_1$.  Similar computations yield
\begin{eqnarray*}
e_1 &\mapsto& \sqrt{2}\cdot\frac{1}{3}\cdot\frac{\sqrt{1}}{1}\cdot \frac{1}{\sqrt{1}} \, e_0, \\
e_2 &\mapsto& \sqrt{3}\cdot \frac{2}{4}\cdot \sqrt{\frac{3}{2}}\cdot \frac{1}{\sqrt{3}}\, e_2, \quad \mbox{\rm and} \\
e_n &\mapsto& \sqrt{n+1} \cdot\frac{n}{n+2} \cdot \sqrt{\frac{2(n-1)+1}{n-1+1}}\cdot \frac{1}{\sqrt{2(n-1)+1}} \, e_{2(n-1)}, \quad n \geq 1.
\end{eqnarray*}
It is productive to simplify this last to
\begin{equation} \label{enmapsto}
e_n \mapsto \sqrt{\frac{n}{n+2}} \cdot \sqrt{\frac{n+1}{n+2}} \, e_{2n-2}, \quad n \geq 3.
\end{equation}

The above makes it clear that the space $\bigvee \{e_0, e_1\}$ is reducing for $B^*_{z^2 \overline{z}}$ with the restriction to the subspace acting via the matrix
$$\left(\begin{array}{cc}
0 & \frac{\sqrt{2}}{3}\\
\frac{\sqrt{2}}{3}& 0
\end{array}\right).$$
Observe that this matrix is self adjoint and so normal.
Likewise, the space $\bigvee\{e_2\}$ is reducing and of course on this space $B^*_{z^2 \overline{z}}$ is normal.

Since the map $f$ defined by $f(n) = 2(n-1)$ is injective on $\mathbb{N}_3$ it is easy to see that $\mathbb{N}_3$ may be partitioned into (disjoint) sets each of which is the orbit $\mathcal{O}_j$ under iterates of $f$ of a ``founder'' $j$:
$$\mathcal{O}_j = \{f^n(j): n \in \mathbb{N}_0 \},$$
where $f^n$ denotes the $n$-fold composition of $f$ with itself.
An easy induction shows that
\begin{equation}  \label{eq:goodformOhsubj}
\mathcal{O}_j = \{2^n(j-2) + 2 : n \in \mathbb{N}_0\}.
\end{equation}
(Given some $m \in \mathbb{N}_3$, to find the founder of the orbit in which it resides, choose the largest $j \in \{0, 1, 2, \ldots \}$ such that with $g(k) = \frac{k}{2} + 1$, $g^j(m)$ is a non-negative integer.)  Let $\mathcal{F}$ denote the collection of integers which are founders.  For example, $3$ is a founder, and $\mathcal{O}_3 =  \{3, 4, 6, 10, \ldots\}$ while $\mathcal{O}_5 =  \{5, 8, 14, 26, \ldots\}$.  Finally, it is clear that there are countably many such (disjoint) orbits, since (at least) every odd number in $\mathbb{N}_3$ must be a founder because the range of $f$ is contained in $2 \mathbb{N}_3$.

If we fix temporarily a founder $j$ in $\mathbb{N}_3$ it is clear that $\mathcal{H}_j := \bigvee \{e_n: n \in \mathcal{O}_j \}$ is reducing for $B^*_{z^2 \overline{z}}$.  Further, on this space $B^*_{z^2 \overline{z}}$ acts as a weighted shift.  Taking for another example $\mathcal{H}_5$, and renumbering the basis as $(h_i)_{i=0}^\infty$ so that $h_i = e_{f^i(5)}$, one calculates that the weights of the shift are $w_n = \sqrt{\frac{2^n + 2}{2^n + 4}} \cdot \sqrt{\frac{2^n + 3}{2^n + 4}}$.

Each of the terms in this product is a GRWS;  the second is in the extension of a portion of Sector IV of the magic square (the wedge bounded by $D = N$ and $D = p N$ with $p = 2$, but now with $N$ and $D$ both larger than $1$).  Were this the only weight, the shift would be hyponormal but not $2$-hyponormal, since while the general discussion in \cite{BCE4} assumes that $-1 < N < 1$ and $-1 < D < 1$, the argument for this sector is based on determinants and does not require this assumption.  The first term is the weights for a subnormal shift on the special line $D = 2 N$ from Sector IV, but extended out of the magic square;  one may check that as suggested by the work for the magic square, the Berger measure for this shift is given by
$$\frac{1}{5} \cdot (4 \delta_{1/2} + \delta_1),$$
where $\delta_a$ is the usual Dirac measure at $a$.
One may also check by examination of the three by three matrices relevant for $2$-hyponormality that the product of these two shifts is not $2$-hyponormal.  Since the Schur product of hyponormal shifts is hyponormal, the product shift is;  since this is part of a direct sum, the direct sum of the shifts cannot be $2$-hyponormal.  It will follow from the general results that each of the shifts in the sum is at least hyponormal, and thus their direct sum is.  Thus the operator $B^*_{z^2 \overline{z}}$ has an orthogonal decomposition of the form
$$B^*_{z^2 \overline{z}} = \mbox{\rm M}_2 \oplus D_1 \oplus \left(\oplus_{j \in \mathbb{N}_1 \cap \mathcal{F}} W(j)\right)$$
where M$_2$ is a two by two normal matrix, $D_1$ is a (real) matrix on $\mathbb{C}$ and hence normal, and the $W(j)$ are the weighted shifts on the $\mathcal{O}_{j}$ all of which are hyponormal and at least one of which is not $2$-hyponormal.

If we fix $j$ some founder for a moment, and set the vectors $u_n:= e_{f^n(j)}$ to be a basis for the acting space, then using \eqref{enmapsto} and \eqref{eq:goodformOhsubj}, we obtain
\begin{eqnarray*}B_{z^2\bar z}^* u_n &=& \sqrt{\frac{ 2^n(j-2)+2}{2^n(j-2)+2+2}}  \sqrt{\frac{2^n(j-2)+2
+1}{2^n(j-2)+2+2}} u_{n+1}\\
&=& \sqrt{\frac{ 2^n+\frac{2}{j-2}}{2^n+ \frac{4}{j-2}  }}  \sqrt{\frac{2^n+  \frac{3}{j-2}}{2^n+ \frac{4}{j-2}}} u_{n+1}.
\end{eqnarray*}
\bigskip
Thus the general expression of the weights for the shift corresponding to $\mathcal{O}_{j}$ is of the form
$$\sqrt{\frac{2^n + 3/(j-2)}{2^n + 4/(j-2)}} \cdot \sqrt{\frac{2^n + 2/(j-2)}{2^n + 4/(j-2)}}, \quad n \geq 0, j \geq 3.$$
Considerable experimentation using \textit{Mathematica} (\cite{Wol}) suggests that in fact none of these (product) shifts are $2$-hyponormal.

Observe finally for this example that each of the shifts is contractive, $D_1$ is contractive, and so is M$_2$, so $B^*_{z^2 \overline{z}}$ is a contractive operator.

This example completed, we turn next to the general results.

\medskip

\subsubsection{$B^*_{z^t \overline{z}^s}$ with $s \geq t$}  \label{subse:zuptzbzrupssgeqt}\hfill\\

Recall that we have
$$B^*_{z^t \overline{z}^s} = K^* P_{harm}(M_{z^s \overline{z}^t}).$$
(Here and in what follows, when the symbol is understood we will sometimes simply write this as $B^*$ and similarly write the adjoint as $B$ to ease the notation.)  In this case $s \geq t$, things are in some sense ``analytic,'' because we will end up applying $P_{harm}$ only to products $z^{m+s} \overline{z}^t$ with $s + m \geq t$ and thus using only the first line of $P_{harm}$'s definition. We will set $d = s-t$, which will help exhibit which portions of the resulting structure results depend on $s$ and which only on the difference $d = s-t \geq 0$.  It will turn out that the form of the direct sum decomposition, although not the ``weights,'' will depend only on $d$, and this will later allow for results considering a sum like $B^*_{z^t \overline{z}^s} + B^*_{z^k \overline{z}^m}$ so long as $s-t = m-k$;  see similar results in \cite{KLL}.

One computes that
\begin{equation*}  \label{eq:Bstarsgtzupm}
B^*_{z^t \overline{z}^s}(z^m) = \left(\frac{d + m + 1}{s + m + 1}\right) \cdot \sqrt{\frac{2(d + m) + 1}{d + m + 1}} z^{2(d + m)}
\end{equation*}
and then, using that the orthonormal basis has $e_m = \sqrt{m+1} z^m$, that
\begin{eqnarray}  \label{sgttmapsto}
B^*_{z^t \overline{z}^s}(e_m) &=& \sqrt{m+1} \left(\frac{d + m + 1}{s + m + 1}\right) \cdot \sqrt{\frac{2(d + m) + 1}{d + m + 1}} \frac{1}{\sqrt{2(d+m) + 1}}e_{2(d + m)} \nonumber \\
&=& \sqrt{\frac{m+1}{s + m + 1}} \cdot \sqrt{\frac{d+m+1}{s + m + 1}} e_{2(d + m)}.   \label{eq:Bstaractonemsgeqt}
\end{eqnarray}

The map $f$ defined by $f(m) = 2(d + m)$ is clearly injective on  $\mathbb{N}_0$, and so $\mathbb{N}_0$ may be partitioned into disjoint orbits $\mathcal{O}_m$ induced by a founder $m$ (such that $m/2 - d$ is not a non-negative integer) and so that $\mathcal{O}_m= \{f^k(m): k \in \mathbb{N}_0\}$.  Observe that there are infinitely many such orbits;  $0$ induces an orbit, as does any odd positive integer $m$ (since the range of $f$ is contained in $2 \mathbb{N}_0$), and depending on the size of $d$ some even $m$ may be founders.

In the special case $d = 0$ the orbit induced by the founder $0$ consists only of $0$, and one obtains that
\begin{equation}  \label{eq:e0mapstoifdeq0}
e_0 \mapsto \frac{1}{s + 1} e_0 \quad \mbox{\rm (if $d$ = 0)}.
\end{equation}
All other orbits when $d = 0$ are infinite, as are all orbits for $d > 0$.

If we consider the basis vectors arising from an orbit $\mathcal{O}_j$ with (fixed) founder $j$, namely $\{e_m: m = f^k(j), k \in \mathbb{N}_0\}$ (and excluding the anomalous case $d = 0$ and $j = 0$), the restriction of $B^*_{z^t \overline{z}^s}$ to the (reducing) subspace of their span acts as a weighted shift as in \eqref{eq:Bstaractonemsgeqt}.    One computes again using an induction, and that $f(m) = 2(d + m)$, that $\mathcal{O}_j:=\{2^n(2d+j)-2d\ : \ n\in \mathbb{N}_0\}$. If we set the vectors $v_n:= e_{g^n(j)}$, then using  \eqref{sgttmapsto} we obtain
\begin{eqnarray*}
B_{z^t\bar z^s}^* v_n &=& \sqrt{\frac{ 2^n(2d+j)-2d+1}{ s+2^n(2d+j)-2d +1    }}  \sqrt{\frac{ d+ 2^n(2d+j)-2d+1       }{s+ 2^n(2d+j)-2d+1}} v_{n+1} \\
 &=& \sqrt{\frac{ 2^n+\frac{1-2d}{2d+j}}{2^n+ \frac{1+s-2d}{2d+j}  }}  \sqrt{\frac{2^n+  \frac{1-d}{2d+j}}{2^n+ \frac{1+s-2d}{2d+j}}} v_{n+1}.
\end{eqnarray*}
Thus the weights of this shift we denote $W_j$ satisfy
\begin{eqnarray}
w_n&=& \sqrt{\frac{2^n(2d + j)+(1-2d)}{2^n(2d + j)+(1+s-2d)}}\cdot \sqrt{\frac{2^n(2d + j)+(1-d)}{2^n(2d + j)+(1+s-2d)}} \nonumber \\
&=&  \sqrt{\frac{2^n+\frac{1-2d}{2d+j}}{2^n+\frac{1+s-2d}{2d+j}}}\cdot \sqrt{\frac{2^n+\frac{1-d}{2d+j}}{2^n+\frac{1+s-2d}{2d+j}}}.  \label{eq:weightsforsgeqt}
\end{eqnarray}

Note for future use that the orbits, and thus the reducing subspaces for the direct sum decomposition, depend upon $d$ and $j$ but not upon $s$.  Thus the decompositions for the same $d$, but different $s$, are conformable and if we consider the addition of various $B^*$ with the same $d$, we must simply examine the one-dimensional normals (if present in the case $d = 0$) and the shifts on each orbit shift space for some property of interest.

These computations and observations in hand, we turn first to the special case $s = t$ and then to the general case $s > t$.

\vspace{.1in}

\noindent \textbf{A special case:} $B^*_{z^s \overline{z}^s}$ ($s = t > 0$, $d = s-t = 0$)

\medskip

In this case one has, as in \eqref{eq:e0mapstoifdeq0},  $e_0 \mapsto \frac{1}{s + 1} e_0$, and so in the direct sum decomposition $e_0$ spans a one-dimensional subspace reducing for $B^*_{z^s \overline{z}^s}$ on which it is normal (this corresponds to the entry $D_1$ in the sum for the motivating example;  there is no first matricial term).  For the shift subspaces, the shift $W_j$ arising from some founder $j \geq 1$ has weights, from \eqref{eq:weightsforsgeqt} with $d=0$,
\begin{equation*}  \label{eq:wtsWjfordeq0}
w_n = \left(\sqrt{\frac{2^n + \frac{1}{j}}{2^n + \frac{s+1}{j}}}\right)^2.
\end{equation*}
Each of the two (in this case, identical) terms in the Schur product is the weight sequence for a GRWS, and the subnormality or level of $k$-hyponormality of these terms is known exactly.   We have $N = N(j)= \frac{1}{j} > 0$ and $D = D(j) = (s+1)N(j)$, and there are two subcases. For the first, if $s$ is such that $s+1 = 2^k$ for some $k$ in $\mathbb{N}$, then the resulting shift is subnormal (this is one of the extensions of the special lines in Sector IV of the magic square;  see the discussion after Theorem \ref{th:detailsofmagicsquare}).  Thus $s$ of the form $s =2^k - 1$, for some $k \in \mathbb{N}$, yields that each term is subnormal, and since the Schur product of subnormal shifts is subnormal, so it this restriction of $B^*_{z^s \overline{z}^s}$.  Observing that this is independent of $j$, we obtain subnormality.

For the second subcase,  if $s$ satisfies $2^k-1 < s < 2^{k+1} - 1$ for some $k \in \mathbb{N}$, then each of these terms is known to be exactly $k$-hyponormal but not $(k+1)$-hyponormal.  Thus their Schur product is at least $k$-hyponormal.  Experiments show that the Schur product may have a higher hyponormality, but we do not know what.

Therefore we have the following.

\begin{thm}  \label{th:seqt}
Let $B_{z^s \overline{z}^s}$ be an H-Toeplitz operator on the Bergman space and suppose $s > 0$.  Then if $s =2^k - 1$ for some $k \in \mathbb{N}$, $B^*_{z^s \overline{z}^s}$ is  the direct sum of a one-dimensional normal operator and an infinite family of subnormal weighted shifts as described above, and therefore subnormal.  If $s$ satisfies $2^k-1 < s < 2^{k+1} - 1$ for some $k \in \mathbb{N}$, then $B^*_{z^s \overline{z}^s}$ is the direct sum of a one-dimensional normal operator and an infinite family of $k$-hyponormal weighted shifts as described above with each entry in the two-term Schur product exactly $k$-hyponormal, and is therefore each shift in the infinite family is at least $k$-hyponormal.  In either case it is clear that $B_{z^s \overline{z}^s}$ is not hyponormal.

Since each of the relevant shifts is contractive, and the entry for the one-dimensional normal operator is $\frac{1}{s+1}$, it is clear that both $B^*_{z^s \overline{z}^s}$ and $B_{z^s \overline{z}^s}$ are contractive.
\end{thm}

Using that the decompositions into reducing subspaces depend only on $d$ and are compatible, we have the following.  Let $\ell^1$ denote the space of absolutely summable complex sequences.

\begin{thm}  \label{th:sumswhendeq0}
Let $(a_n)_{n=0}^\infty$ be a sequence of real, non-negative, numbers in $\ell^1$. Let $\{B_{z^{s_i} \overline{z}^{s_i}}, s_i \in \mathbb{N}, i \in \mathbb{N}\}$ be a collection of H-Toeplitz operators on the Bergman space.  Define $f$ by $f(z, \overline{z}) = \sum_{i=1}^\infty a_i z^{s_i} \overline{z}^{s_i}$.  Then $B^*_{f(z, \overline{z})}$ is hyponormal.  Therefore, $B_{f(z, \overline{z})}$ is not hyponormal.  If $\sum_{i=0}^\infty a_i \leq 1$,  then $B^*_{f(z, \overline{z})}$ and $B_{f(z, \overline{z})}$ are contractive.
\end{thm}

\noindent Proof.  Since hyponormality for a weighted shift is equivalent to non-decreasing weights, it is clear that the sum of hyponormal weighted shifts is hyponormal.  Of course the sum of the normal (one-dimensional) portions of the decomposition is normal.  Using the compatibility of the decompositions, the result follows.  \pfend

It is reasonable to ask whether something strong than hyponormality holds, at least if the individual summands are $k$-hyponormal or even subnormal.  This seems to be difficult;   for example, the sum of two subnormal operators, even two GRWS on the same special line, need not be subnormal (an example is the pair with parameters $(p,N,D) = (2,1/3, 2/3)$ and  $(p,N,D) = (2,1/5, 2/5)$).

\vspace{.1in}

\noindent \textbf{The general case:} $B^*_{z^t \overline{z}^s}$ ($s > 0$, $d = s-t > 0$)

\medskip

In this case there is no one-dimensional reducing subspace in the decomposition, and everything comes down to the weighted shifts.  It turns out that we get subnormality in all cases, but in some cases we get even more, such as $\mathcal{MID}$ shifts.  Since the weights for the shift $W_j$ are a product, namely
\begin{equation}   \label{eq:weightsforsgtanddg1}
w_n =  \sqrt{\frac{2^n+\frac{1-2d}{2d+j}}{2^n+\frac{1+s-2d}{2d+j}}}\cdot \sqrt{\frac{2^n+\frac{1-d}{2d+j}}{2^n+\frac{1+s-2d}{2d+j}}},
\end{equation}
there are various combinations of the ``goodness'' of the two shifts.

For the first term in the product, the parameter ``$N$'' of the GRWS form, which is $\frac{1-2d}{2d+j}$,  is always negative and easily seen to be greater than $-1$ for every $j$, and thus so far we are in the magic square.  If for the choice of $s$ and $d$ we have $1 + s - 2d \leq 0$, then we are in Sector I and the shift is $\mathcal{MID}$ (and in fact, as will be used later, even better in that its weights squared are interpolated by a Bernstein function).  If $1 + s - 2d \geq 0$ but $1 + s - 2d \leq 2d - 1$, i.e., $s \leq 4d - 2$, we are in Sector II and $\mathcal{MID}$.  If  $1 + s - 2d \geq 2d - 1$, so $s > 4d-2 = 2(2d-1)$, we are in Sector III or vertically above Sector III (depending on $j$) and have subnormality (since, as noted above, the measure arguments in \cite{BCE1} for Sector III go through without change for $D \geq 1$ leaving the magic square).  Summarizing, the first shift of the product is
$$\begin{array}{cc}
\mathcal{MID} \mbox{\rm (Bernstein, Sector I)}&  s \leq 2d-1\\
\mathcal{MID} \mbox{\rm (Sector II)}&  2d-1 < s \leq 2(2d-1)= 4d-2\\
\mbox{\rm subnormal (Sector III or above)}& s > 2(2d-1).
\end{array}
$$

For the second shift in the product we must consider separately $d = 1$ and $d > 1$.  In the case $d = 1$,
there are two subcases.  If also $s = 1$, then the second shift reduces to the unweighted shift, so the weights in the first term of the product in \eqref{eq:weightsforsgtanddg1} are all that is present, and in fact we have that the resulting shift is $\mathcal{MID}$ because it is on the boundary of Sectors I and II.  If $d =1$ and $s > 1$, then the second shift is subnormal (on the positive $D$ axis, in or outside of the magic square, and again the measure arguments go through without difficulty).

In the case $d \geq 2$, then the analysis proceeds as for the first term yielding the summary for the second shift of the product  (note that $s > d-1$ is automatic)
$$\begin{array}{cc}
\mathcal{MID} \, \mbox{\rm (Bernstein, Sector I)}&  s \leq 2d-1\\
\mathcal{MID} \, \mbox{\rm (Sector II)}&  2d-1 < s \leq 3d-2\\
\mbox{\rm subnormal (Sector III or above)}& s > 3d-2.
\end{array}
$$

The product of $\mathcal{MID}$ weights is $\mathcal{MID}$;  the product of weights squared interpolated by a Bernstein function with another is $\mathcal{MID}$ but need not have weights squared interpolated by a Bernstein function;
it is unknown what can happen with the product of the weights of an $\mathcal{MID}$ shift with another shift of some less restrictive type.  The notion of $\mathcal{MID}$ has been defined for shifts, but it may easily be extended to any operator $T$ whose (infinite) matrix satisfies the condition that, for every $p > 0$,  if every entry is raised to the $p$-th power the operator is subnormal.  Clearly if the operator is a direct sum of shifts this amounts to each of the shifts being $\mathcal{MID}$.

There are thus numerous combinations of the ``goodness'' of the two pieces, and we leave to the interested reader a finer tabulation of the various combinations and content ourselves with the following.

\begin{thm}
Consider some H-Toeplitz shift $B_{z^t \overline{z}^s}$ with $s > 0$ and $d = s-t > 0$.  If $s \leq 3d-2$ then $B^*_{z^t \overline{z}^s}$ is $\mathcal{MID}$ (that is, it is a direct sum of weighted shifts each of which is $\mathcal{MID}$); note that this includes the special case $d = 1$ and $s = 1$.  For all choices of $s$ and $d$ as in the hypotheses, $B^*_{z^t \overline{z}^s}$ is a direct sum of weighted shifts each subnormal, and is therefore subnormal.  It follows that in either case $B_{z^t \overline{z}^s}$ is not hyponormal.
\end{thm}

We remark that for $d = 1$ and $s = 2$ we have, for each of the shifts in the reducing subspaces for $B_{z^t \overline{z}^s}$, the product of an $\mathcal{MID}$ shift and a ``merely'' subnormal shift, and it is not in general known whether such a shift, surely subnormal, can in fact be $\mathcal{MID}$.

We have the analog of Theorem \ref{th:sumswhendeq0}.

\begin{thm}  \label{th:sumsofHNliked}
Consider some $d > 0$ fixed, and a family of H-Toeplitz shifts $\{B_i := B_{z^{t_i} \overline{z}^{s_i}}, i \in \mathbb{N}\}$ with, for all $i$, $s_i > 0$ and $d = s_i-t_i $.  Suppose $(a_i)_{i=1}^\infty$ is a sequence of non-negative numbers in $\ell^1$.  Define $f$ by $f(z, \overline{z}) = \sum_{i=1}^\infty a_i z^{s_i} \overline{z}^{s_i}$.

\begin{itemize}
\item  If $s_i \leq 2d-1$ for all $i$,  then $B^*_{f(z, \overline{z})}$ is $\mathcal{MID}$ (that is, it is a direct sum of weighted shifts each of which is $\mathcal{MID}$).  (Clearly with fixed $d$ it is equivalent to take the family and the sequence $a$ to be finite.)
\item For all choices of the $s_i$ and $t_i$ as in the hypotheses, $B^*_{f(z, \overline{z})}$ is a direct sum of weighted shifts each hyponormal, and is therefore hyponormal.
    \end{itemize}

It follows in either case that $B_{\sum_{i=1}^\infty a_i z^{t_i} \overline{z}^{s_i}}$ is not hyponormal.
\end{thm}

\noindent Proof.  The second assertion is as in the proof of Theorem \ref{th:sumswhendeq0}.  For the first, if a shift has weights squared completely alternating (in particular, if the weights squared are interpolated by a Bernstein function), then the sequence of weights is completely alternating (the class of completely alternating sequences is closed, for any $0 \leq p \leq 1$,  under the operation of raising each term to the $p$-th power:  see \cite[Corollary 3]{AR}, or see \cite[Chapter 3, Corollary 2.10]{BCR}).  Clearly a (convergent) sum of any number of completely alternating sequences is completely alternating.  But if a sequence is completely alternating, it is log completely alternating.  None of this is changed by multiplying by (positive) constants, and log completely alternating weights are sufficient for $\mathcal{MID}$. \pfend

Finally, observe that each of the shifts in the direct sum decomposition is contractive (even should one of the terms in the product for the weights depart from the magic square), as is the normal direct summand.  Therefore, we have

\begin{prop}
Consider some H-Toeplitz operator $B_{z^t \overline{z}^s}$ with $s > 0$ and $d = s-t > 0$.  Then $B_{z^t \overline{z}^s}^*$, and hence $B_{z^t \overline{z}^s}$, are contractive.  For some $d > 0$ fixed, and a family of H-Toeplitz shifts $\{B_i := B_{z^{t_i} \overline{z}^{s_i}}, i \in \mathbb{N}\}$ with, for all $i$, $s_i > 0$ and $d = s_i-t_i $, and given $(a_i)_{i=1}^\infty$ a sequence of non-negative numbers with sum less than or equal to $1$, define $f$ by $f(z, \overline{z})=\sum_{i=1}^\infty a_i z^{t_i} \overline{z}^{s_i}$. Then $B^*_{f(z, \overline{z})}$  is contractive, as is $B_{f(z, \overline{z})}$.
\end{prop}

\medskip

\subsubsection{$B^*_{z^t \overline{z}^s}$  with $t > s$}\hfill\\

\medskip

We turn next to the case, as in the motivating example beginning Subsection \ref{subse:HToepinBerg}, in which the coanalytic term in the product is of higher degree.  Set $\delta = t-s$ and assume $\delta > 0$.  We assemble some computational results.

\begin{lem}   \label{le:computformorez}
Consider $B^*_{z^t \overline{z}^s}$ with $t > s$ and with $\delta = t-s$.  There are four cases (where we temporarily simplify $B^*_{z^t \overline{z}^s}$ to $B^*$ to ease the notation):
\begin{enumerate}
\item  If $0 \leq m \leq \delta -1$, then $B^*(z^m) = \frac{\delta-m + 1}{s + \delta + 1} \sqrt{\frac{2(\delta-m)}{(\delta-m) + 1}} z^{2(\delta-m) - 1}$, so
    \newline $B^*(e_m) = \sqrt{m+1}\frac{\delta-m + 1}{s + \delta + 1} \sqrt{\frac{2(\delta-m)}{(\delta-m) + 1}} \frac{1}{\sqrt{2(\delta-m)-1+1}}
e_{2(\delta-m)-1} =  \frac{\sqrt{m+1}\sqrt{\delta-m+1}}{\delta+s + 1}e_{2(\delta-m)-1}$;
\item If $\delta   \leq m \leq 2\delta -1$, then $B^*(z^m) = \frac{m-\delta + 1}{m + s + 1} \sqrt{\frac{2(m-\delta)+1}{m-\delta + 1 }} z^{2(m-\delta)}$, so
    \newline $B^*(e_m) = \sqrt{m+1}\frac{m-\delta + 1}{m + s + 1} \sqrt{\frac{2(m-\delta)+1}{m-\delta +1}} \sqrt{\frac{1}{2(m-\delta) + 1}}e_{2(m-\delta)} = \frac{\sqrt{m+1}\sqrt{m-\delta+1}}{ m +s + 1}e_{2(m-\delta)}$;
\item If $m = 2 \delta$, then $B^*(z^m) = B^*(z^{2\delta}) = \frac{\delta + 1}{2 \delta +s + 1}\sqrt{\frac{2 \delta + 1}{\delta + 1}}z^{2 \delta}$, so
      \newline $B^*(e_{2 \delta}) =  \sqrt{2 \delta + 1} \frac{\delta + 1}{2 \delta +s +  1}\sqrt{\frac{2 \delta + 1}{\delta + 1}}\frac{1}{\sqrt{2\delta + 1}} e_{2 \delta} = \frac{\sqrt{2 \delta + 1} \sqrt{\delta + 1}}{2 \delta + s + 1} e_{2 \delta}$;
\item If $m \geq 2 \delta + 1$, then $B^*(z^m) = \frac{m - \delta + 1}{ m + s +1}\sqrt{\frac{2(m-\delta) + 1}{m - \delta + 1}}z^{2(m-\delta)}$, so \newline
          $B^*(e_m) = \sqrt{m+1}\frac{m - \delta + 1}{ m + s +1}\sqrt{\frac{2(m-\delta) + 1}{m - \delta + 1}} \frac{1}{\sqrt{2(m-\delta) + 1}} e_{2(m-\delta)}=\frac{\sqrt{m + 1} \sqrt{m - \delta + 1}}{m + s + 1}e_{2(m-\delta)}$.
\end{enumerate}
\end{lem}

\noindent Proof.  These are simply computations using \eqref{eg:defPharm}, \eqref{eq:Kstarzupn}, and \eqref{eq:Kstarzbarupn} of Lemma \ref{le:computformorez}.  \pfend

One checks that from this result $B^*$ has reducing subspaces $\bigvee \{e_0, \ldots, e_{2 \delta - 1}\}$, $\bigvee \{e_{2 \delta}\}$, and $\bigvee \{e_j : j \geq 2\delta + 1\}$.  From  (1) and (2)   with slightly more work one can see that each row and each column of the matrix for the restriction of $B^*$ to the first space has exactly one non-zero entry and that $B^*B$ and $B B^*$ are diagonal.  Further, one computes for the special case $\delta = 1$ this restriction is not hyponormal for any $s \geq 0$, and further that for $\delta \geq 2$,
$$
\begin{aligned}
	\langle (B^*B-BB^*)e_{2\delta-1},e_{2\delta-1}\rangle&=\|B(e_{2\delta-1})\|^2-\|B^{\ast}(e_{2\delta-1})\|^2\\
	&=\frac{\delta+1}{(s+\delta+1)^2}-\frac{2\delta^2}{(s+2\delta)^2}\\
	&=\frac{(\delta+1)(s+2\delta)^2-2\delta^2(s+\delta+1)^2}{(s+\delta+1)^2(s+2\delta)^2}\\
   &=\frac{-2\delta^4-4s\delta^3-(2s^2-2)\delta^2+(s^2+4s)\delta+s^2}{(s+\delta+1)^2(s+2\delta)^2}\\
   &=\frac{-(\delta-1)(2\delta^3+(4s+2)\delta^2+2s(s+2)\delta+s^2)}{(s+\delta+1)^2(s+2\delta)^2}\\
\end{aligned}
$$
since $B^*(e_{2\delta-1})=\frac{\sqrt{2}\delta}{2\delta+s}e_{2(\delta-1)}$ and
$B(e_{2\delta-1})=\frac{\sqrt{\delta+1}}{s+\delta+1}e_0.$
Since this is obviously negative for any $s \geq 0$, the restriction is never hyponormal for $\delta \geq 1$ and any $s \geq 0$.  One computes directly for the special cases $\delta =1, 2, 3$ that the restriction is not co-hyponormal for any $s$, and as well that for $\delta \geq 4$,
$$
\begin{aligned}
	\langle (B^*B-BB^*)e_{2\delta-2},e_{2\delta-2}\rangle&=\|B(e_{2\delta-2})\|^2-\|B^{\ast}(e_{2\delta-2})\|^2\\
	&=\frac{2\delta^2}{(s+2\delta)^2}-\frac{(2\delta-1)(\delta-1)}{(s+2\delta-1)^2}\\
	&=\frac{2\delta^2(s+2\delta-1)^2-(2\delta-1)(\delta-1)(s+2\delta)^2}{(s+2\delta)^2(s+2\delta-1)^2}\\
	&=\frac{4\delta^3+(8s-2)\delta^2+(3s^2-4s)\delta-s^2}{(s+2\delta)^2(s+2\delta-1)^2},\\
\end{aligned}
$$
since $B^*(e_{2\delta-2})=\frac{\sqrt{2\delta-1}\sqrt{\delta-1}}{s+2\delta-1}e_{2(\delta-2)}$ and
$B(e_{2\delta-2})=\frac{\sqrt{2}\delta}{s+2\delta}e_{2\delta-1}$.

The right-hand side of the commutator expression is clearly positive and so the restriction is never co-hyponormal for $\delta \geq 1$ and any $s \geq 0$.  Thus on this (reducing) subspace the restriction of $B^*$ is neither hyponormal nor co-hyponormal.

From (3) of Lemma \ref{le:computformorez} it is clear that the restriction of $B^*$ to $\bigvee \{e_{2 \delta}\}$ is normal.  As in earlier work, we will consider the restriction of $B^*$ to $\bigvee \{e_j : j \geq 2\delta + 1\}$ and find that it is a direct sum of unilateral shifts, which, since they will be GRWS, we can analyze.

The discussion of orbits, founders, and the like is exactly as in Subsection \ref{subse:zuptzbzrupssgeqt}, \textit{mutatis mutandis}.  One computes that the $n$-th element in the orbit with founder $j$, $j \geq 2 \delta + 1$, is
$$2^n(j - \delta) + 2 \delta.$$
It then follows that the $n$-th weight of the associated shift is
\begin{equation}  \label{eq:weightsforzuptzbarupstgeqs}
\sqrt{\frac{2^n +\frac{2 \delta + 1}{j - 2 \delta}}{2^n + \frac{2 \delta + s + 1}{j - 2 \delta}}}\cdot \sqrt{\frac{2^n +\frac{ \delta + 1}{j - 2 \delta}}{2^n + \frac{2 \delta + s + 1}{j - 2 \delta}}}.
\end{equation}
For the first of these terms, call the associated weighted shift $W_1(s, \delta, j)$;  this is a GRWS with $p = 2$, $N = \frac{2 \delta + 1}{j - 2 \delta}$, and $D = \frac{2 \delta + s + 1}{j - 2 \delta}$, and these fall for any $j \geq 2 \delta + 1$ in Sector IV of the magic square or the extension of that wedge outside the magic square.  Straightforward computations using the results for Sector IV in Theorem \ref{th:detailsofmagicsquare}  yield the following (recalling that the characterization of $k$- but not $(k+1)$-hyponormality extends outside of the magic  square for this sector, since it depends on determinant calculations not restricted to the magic square):

\begin{prop}
For the weighted shift $W_1(s, \delta, j)$ as above, we have that
\begin{enumerate}
\item $W_1(s, \delta, j)$ is subnormal if $s = 0$ or $s = 2^{k-1}(2 \delta + 1)$ for some $k = 1, 2, \ldots$;
\item $W_1(s, \delta, j)$ is $k$-hyponormal but not $(k+1)$-hyponormal if  \newline
\hspace*{.1in} $2^{k-1}(2 \delta + 1) < s < 2^{k}(2 \delta + 1)$, $k = 1, 2, \ldots$.
\end{enumerate}
In particular, $W_1(s, \delta, j)$ is always at least $1$-hyponormal.
\end{prop}

For the second of the terms in \eqref{eq:weightsforzuptzbarupstgeqs} denote the associated shift by $W_2(s, \delta, j)$.  Similar computations yield the following.

\begin{prop}
For the weighted shift $W_2(s, \delta, j)$ as above, we have that
\begin{enumerate}
\item $W_2(s, \delta, j)$ is subnormal if  $s = 2^{k-1}(2\delta-1) + 2^k-1$ for some $k = 1, 2, \ldots$;
\item $W_2(s, \delta, j)$ is $1$-hyponormal but not $2$-hyponormal if $0 \leq s < 1$, and $k$-hyponormal but not $(k+1)$-hyponormal if  \newline
    \hspace*{.1in} $2^{k-2}(2\delta-1) + 2^{k-1}-1 < s < 2^{k-1}(2\delta-1) + 2^k-1$,  $k =  2, 3, \ldots$.
\end{enumerate}
In particular, $W_2(s, \delta, j)$ is always at least $1$-hyponormal.
\end{prop}

One verifies easily that it is not possible for $s$ to satisfy the subnormality conditions for $W_1(s, \delta, j)$ and $W_2(s, \delta, j)$ simultaneously (except for some special cases,  $2^{k-1}(2 \delta + 1)$ is even while $2^{k-1}(2\delta-1) + 2^k-1$ is odd).  It is clear as well that the shift with weights the product in \eqref{eq:weightsforzuptzbarupstgeqs} is always at least $1$-hyponormal and will be at least $k$-hyponormal if $s \geq \min\{2^{k-1}(2 \delta + 1),2^{k-2}(2\delta-1) + 2^{k-1}-1\}$.  We do not know how to analyze in any general way whether/when the product of two $k$-hyponormal shifts might have some hyponormality improved to $(k+1)$ or better, or perhaps even to subnormality (although we believe that if neither shift is subnormal this latter is not possible).  We leave the most detailed analysis of ranges in $s$ for various hyponormality combinations of the two shifts to the reader and content ourselves with the following.

\begin{thm}
Consider $B^*_{z^t \overline{z}^s}$ with $t > s$ and set $\delta = t-s$.  Then $B^*_{z^t \overline{z}^s}$ has a direct sum decomposition as the sum $M \oplus D_1 \oplus \left(\oplus_{j \geq 2 \delta + 1, j \in \mathcal{F}} (W_1(s, \delta, j) \circ W_2(s, \delta, j))\right)$, where $M$ is an operator on $\mathbb{C}^{2 \delta}$ neither hyponormal nor co-hyponormal (except in the case in which $\delta = 1$, when it is normal), $D_1$ is a normal operator on $\mathbb{C}$, and the remainder is a direct sum of weighted shifts each at least $1$-hyponormal.  If $s \geq \min\{2^{k-1}(2 \delta + 1),2^{k-2}(2\delta-1) + 2^{k-1}-1\}$ each of the shifts is at least $k$-hyponormal.  Thus
\begin{itemize}
\item if $\delta = 1$ then $B^*_{z^t \overline{z}^s}$ is at least $1$-hyponormal, and $B_{z^t \overline{z}^s}$ is not hyponormal;
\item  if $\delta > 1$ then neither $B^*_{z^t \overline{z}^s}$ nor $B_{z^t \overline{z}^s}$ is hyponormal.
\end{itemize}

As well, $B^*_{z^t \overline{z}^s}$ and (hence) $B_{z^t \overline{z}^s}$ are contractive.
\end{thm}

\noindent Proof.  One computes directly that if $\delta = 1$ the action of the operator on the two-dimensional space is normal (as in the case of the motivating example).  The only remaining observation is that the shifts and the one-dimensional operator are obviously contractive, and one checks readily that the matrix is contractive.  \pfend

\medskip

We remark that although it is not the point of view of this paper, even when hyponormality (or better) fails for $B^*_{z^t \overline{z}^s}$, since the failure comes from $M$, any such $B^*_{z^t \overline{z}^s}$ is a finite-rank perturbation away from hyponormality.

\medskip

Here (as before) the form of the direct sum decomposition depends on $\delta$ but not on $s$, so that $B^*_{z^{s + \delta} \overline{z}^s}$ and $B^*_{z^{s' + \delta} \overline{z}^{s'}}$ have compatible decompositions.  Using as before that the sum of $1$-hyponormal shifts is $1$-hyponormal the shift portions for an operator like $B^*_{z^{s + \delta} \overline{z}^s + z^{s' + \delta} \overline{z}^{s'}}$ will be hyponormal.  The analysis for hyponormality of such sums then comes down to the hyponormality, or its lack, of the matricial portion of the sum.  While we leave the details to the interested reader, the computations following the proof of Lemma \ref{le:computformorez} show actually that for $\delta \geq 4$ the $(2\delta-1, 2 \delta - 2)$ entry of the matrix for the restriction is larger than the $(0, 2 \delta -1)$ entry, and that the $(2 \delta - 2,2 \delta - 4)$ entry is smaller than the  $(2\delta-1, 2 \delta - 2)$ entry.  Clearly if we compute the weighted sum of such entries, this will remain true, and we will obtain that the matrix in the sum is neither hyponormal nor co-hyponormal.  Thus we have the following.

\begin{thm}
Consider some $\delta > 0$ fixed, and a family of H-Toeplitz shifts $\{B_i := B_{z^{t_i} \overline{z}^{s_i}}, i \in \mathbb{N}\}$ with, for all $i$, $t_i > 0$ and $\delta = t_i - s_i $.  Suppose $(a_i)_{i=1}^\infty$ is a sequence in $\ell^1$ with non-negative terms.  Let $f$ be defined by $f(z, \overline{z}) = \sum_{i=1}^\infty a_i z^{t_i} \overline{z}^{s_i}$.  If $\delta = 1$ the matricial portion of the direct sum decomposition of $B^*_{f(z, \overline{z})}$ is normal, and $B^*_{f(z, \overline{z})}$ is at least $1$-hyponormal.  If $\delta > 1$, then the matricial portion of the direct sum decomposition of $B^*_{f(z, \overline{z})}$ is neither hyponormal nor co-hyponormal, and thus neither $B^*_{f(z, \overline{z})}$ nor $B_{f(z, \overline{z})}$ is hyponormal or co-hyponormal.  As well, if $\sum_{i=1}^\infty a_i \leq 1$, then $B^*_{f(z, \overline{z})}$, and hence $B_{f(z, \overline{z})}$, is contractive.
\end{thm}

\medskip

\subsection{Toeplitz operators on weighted Bergman spaces}\hfill\\

\medskip

In this case we focus on consideration of Toeplitz operators with symbols having a higher analytic degree than co-analytic degree.  (This is because for this space $T^*_{z^p, \overline{z}^q} = T_{z^q, \overline{z}^p}$ in distinction to the H-Toeplitz case considered above.)  As before we assemble some computations for $T_{\overline{z}^sz^{s+d}}(e_n)$ with $s \geq 0$ and $d \geq 0$.

\begin{eqnarray*}
T_{\overline{z}^sz^{s+d}}(e_n)&=&\sqrt{\frac{\Gamma(n+\alpha+2)}{\Gamma(n+1)\Gamma(\alpha+2)}}T_{\overline{z}^sz^{s+d}}(z^n)\cr
&=&\sqrt{\frac{\Gamma(n+\alpha+2)}{\Gamma(n+1)\Gamma(\alpha+2)}}P(\overline{z}^sz^{s+d}z^n)\cr
&=&\sqrt{\frac{\Gamma(n+\alpha+2)}{\Gamma(n+1)\Gamma(\alpha+2)}}\frac{\Gamma(s+d+n+1)\Gamma(s+d+n-s+\alpha+2)}{\Gamma(s+d+n+\alpha+2)\Gamma(s+d+n-s+1)}z^{s+d+n-s}\cr
&=&\sqrt{\frac{\Gamma(n+\alpha+2)}{\Gamma(n+1)\Gamma(\alpha+2)}}\sqrt{\frac{\Gamma(d+n+1)\Gamma(\alpha+2)}{\Gamma(d+n+\alpha+2)}} \, \cdot \\
&& \quad \cdot \, \frac{\Gamma(s+d+n+1)\Gamma(s+d+n-s+\alpha+2)}{\Gamma(s+d+n+\alpha+2)\Gamma(s+d+n-s+1)}e_{d+n}.
\end{eqnarray*}
Put \begin{eqnarray}w(s,d,n)&:=&\sqrt{\frac{\Gamma(n+\alpha+2)}{\Gamma(n+1)\Gamma(\alpha+2)}}\sqrt{\frac{\Gamma(d+n+1)\Gamma(\alpha+2)}{\Gamma(d+n+\alpha+2)}} \times \nonumber\\
&& \times \frac{\Gamma(s+d+n+1)\Gamma(s+d+n-s+\alpha+2)}{\Gamma(s+d+n+\alpha+2)\Gamma(s+d+n-s+1)} \nonumber\\
&=&\sqrt{\frac{(d+n)(d-1+n)(d-2+n)\cdots(1+n)}{(d-1+n+\alpha+2)(d-2+n+\alpha+2)(d-3+n+\alpha+2)\cdots(n+\alpha+2)}}\cr \nonumber \cr
&&\times  \left( \frac{(s-1+d+n+1)(s-2+d+n+1)\cdots(d+n+1)}{(s-1+d+n+\alpha+2)(s-2+d+n+\alpha+2)\cdots(d+n+\alpha+2)}\right)  \nonumber \\
&=&\sqrt{\frac{(d+n)(d-1+n)(d-2+n)\cdots(1+n)}{(d+n+\alpha+1)(d+n+\alpha)(d-1+n+\alpha)\cdots(n+\alpha+2)}}\cr \nonumber \cr
&&\times  \left( \frac{(s+d+n)(s-1+d+n)\cdots(d+n+1)}{(s+d+n+\alpha+1)(s+d+n+\alpha)\cdots(d+n+\alpha+2)}\right).
\end{eqnarray}
\medskip

Of course if $d=0$, this is diagonal (and hence the operator is normal); if $d=1$ this is a weighted shift; if $d>1$ it is a direct sum of $d$ weighted shifts since for $0 \leq i < d$,
$$e_i\mapsto ce_{i+d}\mapsto c^{'}e_{i+2d}\mapsto c^{''}e_{i+3d}\mapsto \cdots.$$
Indeed, for each $0 \leq i < d$ we have an orbit with founder $i$, and on the basis vectors $\{e_{i + \ell d}\}$ the operator acts as a weighted shift.  If in the orbit of some $i$ we consider the weight in moving from $e_{i + \ell d}$ to $e_{i + (\ell+1)d}$, it is $w(s,d,n)$ with $n$ replaced by $i + \ell d$, namely
\small \begin{eqnarray}
&&\sqrt{\frac{(d+i + \ell d)(d-1+i + \ell d)(d-2+i + \ell d)\cdots(1+i + \ell d)}{(d+i + \ell d+\alpha+1)(d+i + \ell d+\alpha)(d-1+i + \ell d+\alpha)\cdots(i + \ell d+\alpha+2)}} \times \nonumber \\
&&\hspace{.2in}\times  \left( \frac{(s+d+i + \ell d)(s-1+d+i + \ell d)\cdots(d+i + \ell d+1)}{(s+d+i + \ell d+\alpha+1)(s+d+i + \ell d+\alpha)\cdots(d+i + \ell d+\alpha+2)}\right). \label{eq:ntonplusd}
\end{eqnarray}\normalsize

But the first of the terms in \eqref{eq:ntonplusd} (the square root) is the product (with order reversed)
\begin{equation*}  \label{eq:firstsqroot}
\prod_{m=1}^{d} \sqrt{\frac{\ell d + m+i}{\ell d + \alpha + m+i + 1}}
\end{equation*}
(Note that $\ell$ is the index relevant to the action of the weighted shift.)  A typical term in this product is, in fact, the weight of a homographic shift $S(A,B,C,D)$ with $A =d$, $B = m+i$, $C = d$, and $D =\alpha + m+i + 1$.    Using that $\alpha > -1$, we see $A D > B C$ and therefore have that each of these shifts is $\mathcal{MID}$ as noted in Subsection \ref {subse:shiftspropertiesclasses}.

The second term in the product in \eqref{eq:ntonplusd} is the Schur square of the product
\begin{equation*}  \label{eq:nonsqrt}
\prod_{j=1}^s \sqrt{\frac{\ell d + j + i + d}{\ell d + j + i + d + \alpha + 1}}.
\end{equation*}
Each of these terms is also the weight for an $S(A,B,C,D)$ with  $A = d$, $B = j + i + d$, $C = d$, and $D = j + i + d + \alpha + 1$;  again this shift is $\mathcal{MID}$.  Finally, since the class of $\mathcal{MID}$ shifts is closed under Schur products of weights, the restriction of $T_{\overline{z}^sz^{s+d}}$ to this shift space is $\mathcal{MID}$ and hence certainly subnormal.  Since $T_{\overline{z}^sz^{s+d}}$ is a direct sum of such shifts, it is subnormal, and $\mathcal{MID}$ in the sense that for any $p > 0$, if we raise every entry in the matricial form for $T_{\overline{z}^sz^{s+d}}$  to the $p$-th power, the result is still subnormal (as a direct sum of subnormal shifts).

Thus we have the following, which improves \cite[Theorem 2.5]{SkL}.

\begin{thm}
Let $T_{\overline{z}^sz^{s+d}}$ be a Toeplitz operator on one of the  weighted Bergman spaces with $s \geq 0$ and $d \geq 0$.  If $d = 0$ the operator is diagonal and normal.  If $d > 1$, the operator is a direct sum of $d$ $\mathcal{MID}$ and hence subnormal weighted shifts;  it is itself $\mathcal{MID}$ in the sense described above, subnormal, and hyponormal. The operator $T_{\overline{z}^sz^{s+d}}$ is contractive.  Except in the case $d = 0$, $T^*_{\overline{z}^sz^{s+d}}$ is not hyponormal.
\end{thm}

\noindent Proof.  It remains only to note that each of the homographic shifts in the decomposition is contractive, and therefore the operator is.  \pfend

\medskip

Observe that the conclusions of subnormality or $\mathcal{MID}$ are independent of the weight $\alpha$ selecting the weighted Bergman space, although the particular weights do reflect $\alpha$.

\medskip

We note that as before the form of the direct sum decomposition depends only on $d$ (a direct sum of $d$ weighted shifts) although the weights depend on $s$.  One obtains the analog of Theorem \ref{th:sumsofHNliked}.

\begin{thm}  \label{thm:sumsforwtdBerg}
Consider some $d > 0$ fixed, and a family of Toeplitz operators  $\{T_i :=T_{\overline{z}^{s_i}z^{s_i+d}}, i \in \mathbb{N}\}$ on a weighted Bergman space with, for all $i$, $s_i \geq 0$ and $d \geq 0$.  Suppose $(a_i)_{i=1}^\infty$ is a sequence of non-negative numbers in $\ell^1$.  Set $T = \sum_{i = 1}^\infty a_i T_i$.  Then $T$ is hyponormal (so $T^*$ is not);  if $\sum_{i=1}^\infty a_i \leq 1$ then $T$ and $T^*$ are contractive.
\end{thm}

From this we can obtain, for example, \cite[Theorem 2.7]{SkL} and (since the argument for \cite[Corollary 2.6]{SkL} is really an argument about hyponormality of a sum of hyponormal $T$ and adjoint $T^*$) we obtain \cite[Corollary 2.6]{SkL} as well.

It is worth noting that the Schur products in \eqref{eq:ntonplusd} are somewhat surprising;  If $\alpha = 0$ and $d = 1$, they turn out to be the Schur product of the Bergman shift weights with the square of the weights of the restriction of the Bergman shift to its first canonical invariant subspace.    Of course one is at liberty to take Schur products of (say) $\mathcal{MID}$ weights, but we were surprised to find this particular Schur product occurring ``in nature.''

\medskip

\subsection{Toeplitz operators on the generalized derivative Hardy space}\hfill\\

In \cite[Theorem 3.2]{KLL} it was shown that, for any non-negative integers $t$ and $s$,  $T_{z^t \overline{z}^s}$ on $S^2_{\alpha, \beta}(\mathbb{D})$ is hyponormal if and only if it is normal if and only if $s = t$ (in which case the operator is easily shown to be diagonal).  We will recapture this result below, via a structure result like those in the previous sections, but clearly we must move away from hyponormality as the question of interest.  We note to the reader that we will enter the world of complete hyperexpansivity and $m$-alternating hyperexpansivity (see Subsection \ref{subse:shiftspropertiesclasses}), and this world is considerably less tractable, and less explored, than the world of subnormality and $k$-hyponormality.  The Schur products of $\mathcal{MID}$ (respectively, subnormal or $k$-hyponormal) shifts are again $\mathcal{MID}$ (respectively, subnormal or $k$-hyponormal); the analogous results in the expansivity world fail.  (For example, the canonical example of a completely hyperexpansive shift is the Dirichlet shift with weights $\sqrt{\frac{n+2}{n+1}}$ -- it is even a $2$-isometry -- yet the shift with weights the (Schur) square of these is not completely hyperexpansive.)   Thus the results in this section will be considerably less ``definitive'' than those in the previous sections.

Familiar computations using Lemma \ref{eq:PforgenderivHS} and \eqref{eq:esubnforgenderivHS}, and with the assumption that $d \geq 1$,  yield that $T_{z^{t + d} \overline{z}^t}$ is the direct sum of $d$ weighted shifts $W_i$, ($i = 0, \ldots, d-1$), where $W_i$ has the weights $w^{(i)}_n$ given by $w^{(i)}_n= w^{(i, \alpha)}_n \cdot w^{(i, \beta)}_n$, where
\begin{equation} \label{eq:wsupisubn}
w^{(i, \alpha)}_n = \sqrt{\frac{d n + i + t + \alpha + d}{d n + i +  \alpha + d}}\sqrt{\frac{d n + i + t + \alpha + d}{d n + i + \alpha}}, \quad n \in \mathbb{N}_0,
\end{equation}
and $w^{(i, \beta)}_n$ is defined similarly.

\medskip

We have the following.
\begin{thm}
Let $T := T_{z^{t + d} \overline{z}^t}$ be a Toeplitz operator on $S^2_{\alpha, \beta}(\mathbb{D})$ where $t \geq 0$ and $d \geq 1$.  Then $T$ is expansive and not hyponormal;  $T^*$ is not hyponormal.
\end{thm}

\noindent Proof.  It is easy to check that the $w^{(i, \alpha)}_n$ as in \eqref{eq:wsupisubn} are both larger than $1$ and decreasing in $n$, as are their $\beta$-analogs $w^{(i, \beta)}_n$;  therefore $T$ is not hyponormal.  Clearly $T^*$, as a direct sum of backward shifts, is not hyponormal.  \pfend

If $d = 1$ we may say a little more, but first must recall a result from \cite{SA}.
\begin{thm}  \label{th:SA}  \cite[pg. 51]{SA}
If a weighted shift has weights $\sqrt{\frac{n + \lambda}{n + 1}}$, then it is $m$-alternatingly hyperexpansive if $m \leq \lambda < m + 1$ and completely hyperexpansive if $1 \leq \lambda \leq 2$.
\end{thm}

\begin{thm}
Let $T := T_{z^{t + 1} \overline{z}^t}$ be a Toeplitz operator on $S^2_{\alpha, \beta}(\mathbb{D})$ where $t \geq 0$.  Then $T$ is a direct sum of  $d$ weighted shifts $W_i$, ($i = 0, \ldots, d-1$), where the weight $w^{(i)}_n$ of $W_i$ is the product $w^{(i, \alpha)}_n \cdot w^{(i, \beta)}_n$;  $w^{(i, \alpha)}_n$ is a two-fold product  $A B$, where $A$ is the weight of a $t$-alternatingly hyperexpansive weighted shift and $B$ is the weight of a $(t + 1)$-alternatingly hyperexpansive weighted shift, and $w^{(i, \beta)}_n$ is a two-fold product  $C D$, where $C$ is the weight of a $t$-alternatingly hyperexpansive weighted shift and $D$ is the weight of a $(t + 1)$-alternatingly hyperexpansive weighted shift.
\end{thm}

Proof.  This is just a computation from Theorem \ref{th:SA} and \eqref{eq:wsupisubn}, and the fact that the restriction of a $m$-alternatingly hyperexpansive shift to one of its canonical invariant subspaces obtained by removing $\bigvee \{e_0, \ldots, e_j\}$ remains $m$-alternatingly hyperexpansive. \pfend

There is more that can be said by identifying the shifts in the decomposition as the reciprocals of well-known shifts.  Each of the weights in the relevant products is the reciprocal of the weights for a subnormal -- in fact, $\mathcal{MID}$ -- homographic shift (see Section \ref{subse:shiftspropertiesclasses}).  Unfortunately, the reciprocal of a subnormal shift need not be completely hyperexpansive: the standard example, due to \cite{At}, is the third Agler shift with weights $\sqrt{\frac{n+1}{n+3}}$, and since this shift is $\mathcal{MID}$ even that is not enough.  We do not know what operator theoretic property being ``reciprocal $\mathcal{MID}$'' provides, so this observation remains somewhat unsatisfying.

\section{Remarks}  \label{se:remarks}

With the structure theorems as above, and in particular the observation that the structure (although not the weights) of the shifts comprising $B^*_{z^m \overline{z}^n}$ depend on $d = m-n$ and not on $m$ and $n$ individually, we may simplify the proofs of various prior results.  We first require two straightforward lemmas.

\begin{lem}  \label{le:scalarmultipofoperatorandHN}
Suppose $T$ is any operator.  Then $T$ is hyponormal if and only if for any constant of modulus $1$ (say, $e^{\mathrm{i} \theta}$), $\hat{T} := e^{\mathrm{i} \theta} T$ is hyponormal.
\end{lem}

Proof. Clearly $\hat{T}^* \hat{T} -\hat{T} \hat{T}^* = T^* T - T T^*$.  \pfend

\begin{lem}
Suppose $W_1$ and $W_2$ are weighted shifts with positive weight sequences $(w^{(1)}_n)_{n=1}^\infty$ and $(w^{(2)}_n)_{n=1}^\infty$ respectively, and suppose in addition each is hyponormal (equivalently, their weight sequences are non-decreasing).  Then for $a, b \in \mathbb{C}$, if $\mathrm{Re}(a \overline{b}) \geq 0$ then $a W_1 + b W_2$ is hyponormal.
\end{lem}

\smallskip

\noindent Proof.  Using the lemma above we may assume $a$ is real and positive;  write $b = x + \mathrm{i} y$ in its usual rectangular form with $x$ and $y$ real.  Clearly $a W_1 + b W_2$ is a unilateral weighted shift;  it is known that any such shift is unitarily equivalent to the shift with weights the moduli of those in the original (\cite[Corollary 1]{Sh}).  It is well-known that hyponormality for a weighted shift is equivalent to the weight sequence non-decreasing.  Putting these together, and using the squares of the moduli for convenience, it suffices to show that
$$|a w^{(1)}_n + b w^{(2)}_n|^2 \leq |a w^{(1)}_{n+1} + b w^{(2)}_{n+1}|^2, \quad n \in \mathbb{N}_0.$$
For a particular $n$, this is
$$(a w^{(1)}_n + x w^{(2)}_n)^2 + y^2 (w^{(2)}_n)^2 \leq (a w^{(1)}_{n+1} + x w^{(2)}_{n+1})^2 + y^2 (w^{(2)}_{n+1})^2,$$
in turn equivalent to
\begin{eqnarray*}a^2 (w^{(1)}_n)^2 &+& x^2 (w^{(2)}_n)^2 + 2 a x w^{(1)}_n w^{(2)}_n + y^2 (w^{(2)}_n)^2 \leq  \\
&\leq& a^2 (w^{(1)}_{n+1})^2 + x^2 (w^{(2)}_{n+1})^2 + 2 a x w^{(1)}_{n+1} w^{(2)}_{n+1} + y^2 (w^{(2)}_{n+1})^2.
\end{eqnarray*}
But using the original weight sequences non-decreasing, and that $a x = \mathrm{Re}(a \overline{b}) \geq 0$, this is obvious.  \pfend

Using the above, we may recapture more transparently the results of \cite[Theorem 2.7]{SkL1} concerning certain sums of analytic Toeplitz operators on the weighted Bergman spaces since matters reduce to sums of shifts.  Theorem 2.8 of that same work considers the hyponormality of the sum of an analytic and a co-analytic Toeplitz operator in that same weighted Bergman setting;  in our notation, what is considered is $a B_{z^m \overline{z}^n} + b B^*_{z^t \overline{z}^s}$ in the case in which $m - n = t - s \geq0$.  The shifts have matching values of $d = m-n = t-s$, and thus the question reduces to (sums of) weighted shifts.  Clearly a necessary condition for hyponormality is non-negativity of the diagonal entries of the matrix of the commutator, and one obtains the result from this point of view somewhat more directly.  One may similarly revisit the situation of H-Toeplitz operators and Theorems 2.9 and 2.10 of \cite{GS}.  We leave such considerations to the interested reader.

Finally, observe that Theorem \ref{th:seqt} includes some shifts which are the Schur product of two shifts, each $k$- but not $(k+1)$-hyponormal.  These are, of course, at least $k$-hyponormal;  experiments suggest that these products may  of higher hyponormality;  we conjecture that, at least, they cannot be subnormal.

\vspace{.1in}

\vspace{.1in}

\noindent {\bf Conflicts of interest}. \  The authors declare no potential conflicts of interest or competing interests.

\vspace{.1in}

\noindent {\bf Acknowledgments}. \  Portions of this work were accomplished during, and materially aided by, visits to Sejong University and Sungkyunkwan University, and we thank the universities and particularly the respective departments for their kind hospitality. \ Several examples in this paper were explored using calculations with the software tool \textit{Mathematica} \cite{Wol}.  The first named author was supported in part by the Labex CEMPI  (ANR-11-LABX-0007-01).  The third named author was supported by the National Research Foundation of Korea (NRF) grant funded by the Korea government (MSIT) (No. 2022R1H1A2091052). The fourth named author was supported by the National Research Foundation of Korea (NRF) grant funded by the Korea government (MSIT)(No. 2021R1C1C1008713).

\end{document}